\documentclass[11pt]{amsart}
\usepackage{palatino}
\usepackage{amssymb,latexsym, amsmath, amsxtra}
\usepackage{graphicx,url}
\usepackage{tikz}
\usepackage[marginratio=1:1,height=8in,width=6.0in]{geometry}

\newtheorem{conjecture}{Conjecture}

\newtheorem{lemma}{Lemma}

\newtheorem{proposition}{Proposition}

\newcommand{\ra}{\rightarrow}
\newcommand{\R}{\mathbf{R}}
\newcommand{\sgn}{\mathrm{sign}}
\begin{document}

\title{Intransitive Dice}

\begin{abstract}
We consider $n$-sided dice whose face values lie between $1$ and $n$ and whose faces sum to $n(n+1)/2$. For two dice $A$ and $B$, define $A \succ B$ if it is more likely for $A$ to show a higher face than $B$. Suppose $k$ such dice $A_1,\dots,A_k$ are randomly selected. We conjecture that the probability of ties goes to 0 as $n$ grows. We conjecture and provide some supporting evidence that---contrary to intuition---each of the $2^{k \choose 2}$ assignments of $\succ$ or $\prec$ to each pair is equally likely asymptotically. For a specific example, suppose we randomly select $k$ dice  $A_1,\dots,A_k$ and observe that $A_1 \succ A_2 \succ \ldots \succ A_k$. Then our conjecture asserts that the outcomes $A_k \succ A_1$ and $A_1 \prec A_k$ both have probability approaching $1/2$ as $n \rightarrow \infty$.
\end{abstract}

\author {Brian Conrey, James Gabbard, Katie Grant, Andrew Liu, Kent E. Morrison}
\address{American Institute of Mathematics,
600 East Brokaw Road, San Jose, CA 95112 USA}

\thanks{
\emph{Mathematics Magazine} \textbf{89} (2016) 133--143.  Research supported in part by a grant from the National Science
Foundation. } 

\maketitle
\renewcommand{\baselinestretch}{1.2}  \large \normalsize 

In 1970 Martin Gardner introduced intransitive (also called ``nontransitive")
dice in his \emph{Mathematical Games} column \cite{Gardner70}. The particular dice he described were invented by Bradley Efron a few years earlier. The six face values for the four dice $A, B, C$, and $D$ are 
\begin{align*}
	A=(0,0,4,4,4,4) \quad &B=(3,3,3,3,3,3) \\
	C=(2,2,2,2,6,6) \quad &D=(1,1,1,5,5,5)
\end{align*}
 The result is a paradoxical cycle of dominance in which 
\begin{itemize}
\item $A$ beats $B$ with probability $2/3$
\item $B$ beats $C$ with probability $2/3$
\item $C$ beats $D$ with probability $2/3$
\item $D$ beats $A$ with probability $2/3$
\end{itemize}
For example, consider rolling $C$ and $D$. Of the $36$ outcomes there are $24$ for which the value shown by $C$ is greater than the value of $D$. 

Effron's dice provide a concrete example of what was first noticed in 1959 by Steinhaus and Trybu{\l}a \cite{SteinhausTrybula59} in a short note (with no mention of dice) showing the existence of independent random variables $X$, $Y$, and $Z$ such that $P(X>Y) > 1/2$, $P(Y>Z) > 1/2$, and $P(Z>X) > 1/2$. This was followed with expanded versions by Trybu{\l}a  containing the details and proofs \cite{Trybula61, Trybula65}. Notably he found equations that describe the maximal probabilities possible for an intransitive cycle of $m$ random variables. For $m=4$ this maximal probability is $2/3$, which means that Efron's dice are optimal in this sense.

Starting with six-sided dice and then generalizing to $n$-sided dice, we focus in this article on just how prevalent intransitive dice are. Much of the work is experimental in nature, but it leads to some tantalizing conjectures about the probability that a random set of $k$ dice, $k \ge 3$, makes an intransitive cycle as the number of sides goes to infinity. For a very restricted ensemble of $n$-sided dice, which we call ``one-step dice,'' we prove the conjectures for three dice.

\section*{How rare are intransitive dice?}
Both surprise and puzzlement are the universal reactions to learning about intransitive dice, and, indeed, that was the case for all of us, but once we had seen some examples, we began to wonder just how special they are.  For example, suppose we pick three dice randomly and find that $A$ beats $B$ and $B$ beats $C$. Does that make it more likely that $C$ beats $A$?

Let's specify exactly what we mean by a random choice of dice. We begin with dice that are much like the standard die commonly used: the number of sides is six, the numbers on the faces come from $\{1,2,3,4,5,6\}$, and the total is $21$. 
We don't  care how the six numbers are placed on the faces and so each die can be represented by a non-decreasing sequence $(a_1,a_2,a_3,a_4,a_5,a_6)$ of integers. Except for the standard die there must be some repetition of the numbers on the faces.  There are 32 such sequences. 

\begin{gather*}
(1, 1, 1, 6, 6, 6), (1, 1, 2, 5, 6, 6), (1, 1, 3, 4, 6, 6), (1, 2, 2,
   4, 6, 6), \\(1, 2, 3, 3, 6, 6), (2, 2, 2, 3, 6, 6), (1, 1, 3, 5, 5, 
  6), (1, 2, 2, 5, 5, 6), \\(1, 1, 4, 4, 5, 6), (1, 2, 3, 4, 5, 6), (2, 
  2, 2, 4, 5, 6), (1, 3, 3, 3, 5, 6),\\ (2, 2, 3, 3, 5, 6), (1, 2, 4, 4,
   4, 6), (1, 3, 3, 4, 4, 6), (2, 2, 3, 4, 4, 6),\\ (2, 3, 3, 3, 4, 
  6), (3, 3, 3, 3, 3, 6), (1, 1, 4, 5, 5, 5), (1, 2, 3, 5, 5, 5),\\(2, 
  2, 2, 5, 5, 5), (1, 2, 4, 4, 5, 5), (1, 3, 3, 4, 5, 5), (2, 2, 3, 4,
   5, 5),\\ (2, 3, 3, 3, 5, 5), (1, 3, 4, 4, 4, 5), (2, 2, 4, 4, 4, 
  5), (2, 3, 3, 4, 4, 5),\\ (3, 3, 3, 3, 4, 5), (1, 4, 4, 4, 4, 4), (2, 
  3, 4, 4, 4, 4), (3, 3, 3, 4, 4, 4)
 \end{gather*}

We say that $A$ \textbf{beats} $B$, denoted by $A \succ B$, if $P(A > B) > P(B > A)$, i.e., the probability that $A > B$ is greater than the probability that $B > A$. Here  we think of the dice as random variables with each of the components in their vector representations being equally likely. This is equivalent to saying that
\[  \sum_{i,j}  \sgn(a_{i}-b_{j}) > 0 .\]
We also say that $A$ dominates $B$ or that $A$ is stronger than $B$. If it happens that $P(A>B) = P(B>A)$, we say that $A$ and $B$ tie or that they have equal strength.

For all choices of three dice $(A,B,C)$ there are $32^3=32{,}768$ possibilities. With the aid of a computer program we found $4417$ triples such that $A \succ B$ and $B \succ C$. Then for the comparison between $A$ and $C$ there were $930$ ties, $1756$ wins for $A$ and $1731$ wins for $C$. Therefore, knowing that $A \succ B$ and $B \succ C$ gives almost no information about the relative strengths of $A$ and $C$. The events $A \succ C$ and $C \succ A$ are almost equally likely!

For each triple of dice there are three pairwise comparisons to make, and for each comparison there are three possible results: win, loss, tie. Throwing away the triples that have any ties leaves us with $13{,}898$ triples and only eight comparison patterns. Our results show that each of the eight patterns occur with nearly the same frequency. Each of the six patterns that give a transitive triple occurs $1756$ times, and each of the two patterns resulting in an intransitive triple occurs $1731$ times. These total $10{,}386$ transitive and  
$3512$ intransitive. 

Rather than look at all ordered triples $(A,B,C)$ we get equivalent information from all subsets of three dice $\{A,B,C\}$, and this, of course, requires much less computation. Of the ${32 \choose 3}=4960$ subsets, we find that $2627$ of them contain ties. Of the remaining $2333$ subsets there are $1756$ transitive sets and $577$ intransitive sets. Allowing for the six permutations of each set we get the same totals as for the ordered triples.

\section*{Proper dice}
With only six sides the number of ties is significant, but what if we increase the number of sides on the dice and let that number grow? Define an $n$-sided die to be an $n$-tuple $(a_1,\dots,a_n)$ of non-decreasing positive integers, $a_1\le a_2\le \dots \le a_n$. 
The {\bf standard} $n$-sided die is 
$(1,2,3,\dots, n)$.
We define {\bf proper} $n$-sided dice to be those with $1 \leq a_i \leq n$ and $\sum a_i = n(n+1)/2$.
Thus, every proper die which is not the standard one has faces with repeated numbers.

Above we listed the proper $n$-sided dice for $n=6$. Here is a list for $n\le5$: 
\begin{eqnarray*}&&
n=1\qquad \qquad (1) \\
&&n=2 \qquad \qquad (1,2)\\
&& n=3 \qquad \qquad 
  (1,2,3),  (2,2,2)
\\
&& n=4 \qquad \qquad (1,1,4,4),  (1,2,3,4),  (1,3,3,3),  (2,2,2,4),  (2,2,3,3)
\\
&&  n=5 \qquad \qquad 
 (1,1,3,5,5),  (1,1,4,4,5),   (1,2,2,5,5),   (1,2,3,4,5),   (1,2,4,4,4), (1,3,3,3,5), \\
&&\qquad  \qquad \qquad \quad (1,3,3,4,4),  (2,2,2,4,5),   (2,2,3,3,5),  (2,2,3,4,4),  (2,3,3,3,4),  (3,3,3,3,3)
\end{eqnarray*}

The number of proper $n$-sided dice occurs as sequence A076822 in the Online Encyclopedia of Integer Sequences \cite{OEIS-1}, where it is described as the number of partitions of the $n$-th triangular number $n(n+1)/2$ into exactly $n$ parts, each part not exceeding $n$. Below are the terms of the sequence for $n \le 27$. 

\[ \begin{array}{rr}
1&1\\ 2& 1 \\3& 2\\4& 5\\5& 12\\ 6& 32\\7& 94\\ 8& 289\\
9& 910\\ 10& 2934\\
11& 9686\\
12& 32540\\
13& 110780\\
14& 381676\\ 
\end{array} 
\qquad \qquad
\begin{array}{rr}
15& 1328980\\ 16& 4669367\\ 17& 16535154\\
18& 58965214\\ 19& 211591218\\
20& 763535450\\
21 &    2769176514\\
22& 10089240974\\
23& 36912710568\\
24& 135565151486\\
25& 499619269774\\
26&  1847267563742\\
27& 6850369296298 \\
&  \\
\end{array}
\]

Obviously, the number of proper dice grows rapidly, and while it is not necessary to our understanding of intransitive dice, we were curious about the rate of growth. Surprisingly, the OEIS entry has nothing about the asymptotics of these partition numbers, but with some heuristics involving the Central Limit Theorem we were able to conjecture that the $n$-th term is asymptotic to \[ \frac{\sqrt{3}}{2\pi}\frac{4^n}{n^2}.\] Eventually we found this result proved rigorously in a 1986 paper by Tak\'acs \cite{Takacs86}, although it is no trivial task to make the connection. You can see the dominant power of 4 in the numbers above. A question that we have been unable to answer is whether there is some construction that gives (approximately) four proper dice with $n+1$ sides from each proper die with $n$ sides.

\section*{Two conjectures for three dice}
Arising from our computer simulations are two conjectures about random sets of three $n$-sided dice as $n \ra \infty$.
\begin{conjecture}
In the limit the probability of any ties is $0$.
\end{conjecture}

\begin{conjecture}
In the limit the probability of an intransitive set is $1/4$.
\end{conjecture}

We can state Conjecture 2 using random ordered triples rather than random sets. For three dice $A, B,C$, there are eight different dominance patterns when there are no ties. In the limit as $n \ra \infty$, we conjecture that each of these patterns has a probability approaching $1/8$. Since two of the patterns are intransitive and six are transitive, the intransitive probability approaches $1/4$ and the transitive probability approaches $3/4$.

Although the conjecture deals with the behavior as $n$ grows, the data for small $n$ already shows us something. For $n=4$, among the ten sets of three distinct dice the only intransitive  set is the set $\{(1,1,4,4), (1,3,3,3), (2,2,2,4)\}$. There is also just one transitive set, while the other eight sets have ties. Thus, the proportion of intransitive is $1/10$. For $n=5$ there are ${12 \choose 3} =220$ sets, and $23$ of these are intransitive with a ratio of $23/220 \approx 0.105$. (There are $54$ transitive sets.) For $n=7$ the computer calculations showed the proportion of intransitive among all the sets is $19929/134044 \approx 0.149$.

A proof of Conjecture 2 appears to be difficult and we do not know how to attack it, but in a later section we present rigorous results about a much smaller collection of dice, the ``one-step dice.``  Conjecture 1, on the other hand, appears to be more attainable, and here we provide a plausibility argument for it. 

First, Conjecture 1, which involves three dice, holds if we can prove that the probability of a tie between two random $n$-sided dice goes to 0 as $n \ra \infty$. That is simply because the probability of any tie is less than or equal to three times the probability of a tie between two dice. We represent a proper $n$-sided die by a vector $(v_1,\ldots,v_n) $ where $v_i$ is the number of faces with the value $i$. There are two restrictions: $\sum v_i = n$, which means that there are $n$ faces, and $\sum i v_i = n(n+1)/2$. Letting $D_n$ be the set of proper $n$-sided dice, we see that $D_n$ is the set of integer lattice points in the intersection of $[1,n]^{n}$ with the $(n-2)$-dimensional affine subspace of $\R^{n}$ defined by the two restrictions. When we roll dice $v$ and $w$, there are $n^{2}$ possible outcomes. When $v$ has the value $i$ showing, then it is greater than $w_{1} + \cdots + w_{i-1}$ of the faces of $w$, and it is less than $w_{i+1}+\cdots +w_{n}$ faces. Summing over $i$, we see that $v$ and $w$ are equally strong when they satisfy the polynomial equation 
\[  \sum_{i < j} v_i w_j  - \sum_{i > j} v_i w_j =0 .\]
This equation defines a hypersurface in $\R^{n} \times \R^{n}$, and 
the set of pairs of tied dice $(v,w)$ is the intersection of  $D_n \times D_n$ with this hypersurface. Heuristically, this means that the dimension of the set of tied dice is one less than the dimension of the set of all pairs. Since the coordinates need to be integers between 0 and $n$,  this suggests that the ratio of the number of tied pairs to the number of all pairs, which is the probability of a tie, should be approximately $1/n$. Computer simulations with $10{,}000$ sample pairs for each size show roughly that behavior.

\[ \begin{array}{ccc}
   n & 1/n & \mbox{Ties} \\
   \hline
   10 &  .100& .0864 \\
   20 &  .050 &  .0329 \\
   30 &  .033 & .0190 \\
   40 &  .025 &  .0124 \\
   50 &  .020 &  .0131 \\
   60 &  .016	 & .0101 \\
   70 &   .014   & .0078 \\
   80  &  .013 & .0061 \\
   90 &   .011  & .0053 \\
   100 &  .010  & .0053
   \end{array}
 \]

A serious difficulty in making this heuristic approach more rigorous is the fact that the both the coordinate range $[1,n]$ and the dimensions of the geometric objects (the affine subspace and the hypersurface) are growing with $n$. However, recent work by Cooley, Ella, Follett, Gilson, and Traldi  \cite{Traldi14} proves that the proportion of ties goes to $0$ as $n \ra \infty$ for $n$-sided dice with values between $1$ and a fixed integer $k$ and having a total equal to $n(k+1)/2$. 

For $n > 7$ we have estimated the probability of intransitive triples by sampling from the sets of proper dice.
Our data, based on 10,000 sample triples  
for each of 
10-sided, 20-sided, 30-sided, 40-sided, and 50-sided
dice, is below. It shows the proportion with ties and the proportion that are intransitive.
\[ \begin{array}{crc}
n & \mbox{ Ties} & \mbox{Intransitive}\\
\hline
10 & .2219 & .1933 \\
20 & .0862 & .2267  \\
30 & .0557 & .2357  \\
40 & .0439 & .2380 \\
50 & .0306 & .2448 
\end{array}
\]

\section*{Other ensembles}
\begin{figure}[h]
\hfill \includegraphics[width=2.5in]{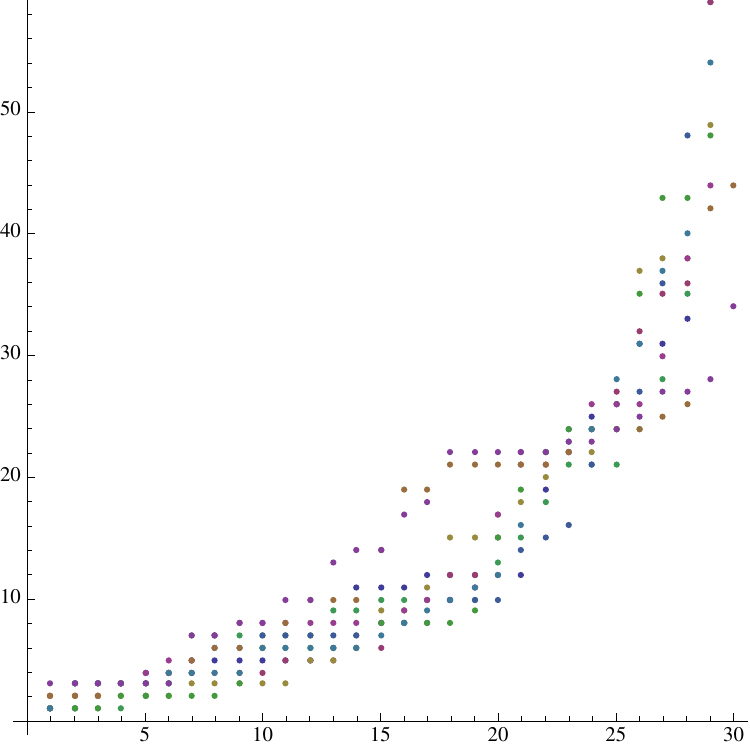}
\hfill \includegraphics[width=2.5in]{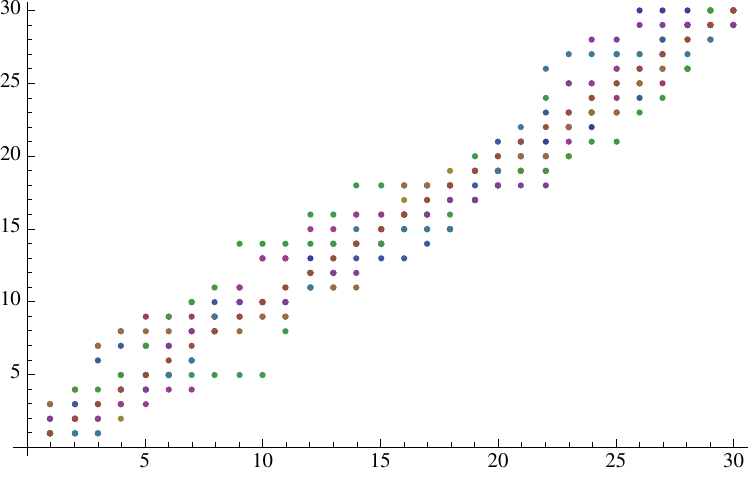}\hspace*{\fill}
\caption{Left: ten random improper dice. Right: ten random proper dice.}
\label{sample-plots}
\end{figure}

We have investigated other ensembles of dice in an effort to see whether the $1/4$ probability of being intransitive is a widespread phenomenon. Suppose we consider $n$-sided dice with the only restriction being that the total is $n(n+1)/2$, thus allowing values greater than $n$. Let's call them ``improper dice.'' These dice are the partitions of $n(n+1)/2$ with exactly $n$ parts. There are significantly more of them. For $n=10$ there are $33{,}401$ improper dice compared with $2934$ proper dice. In a random sample of $1000$ triples of these dice with $20$ sides, we found  $13$ intransitive sets, $958$ transitive sets, and $29$ with one or more ties. These results are far different than for proper dice! What's the cause? You can visualize an $n$-sided die by looking at the plot in the plane of the point set  $\{(i, a_i)\}$. The left side of Figure \ref{sample-plots}  shows the superimposed plots for ten random improper dice with $30$ sides. The right side shows the same for ten random proper dice. You can see that the typical proper die is much closer to the standard die than the typical improper die.

\begin{figure}

\includegraphics[width=2.5in]{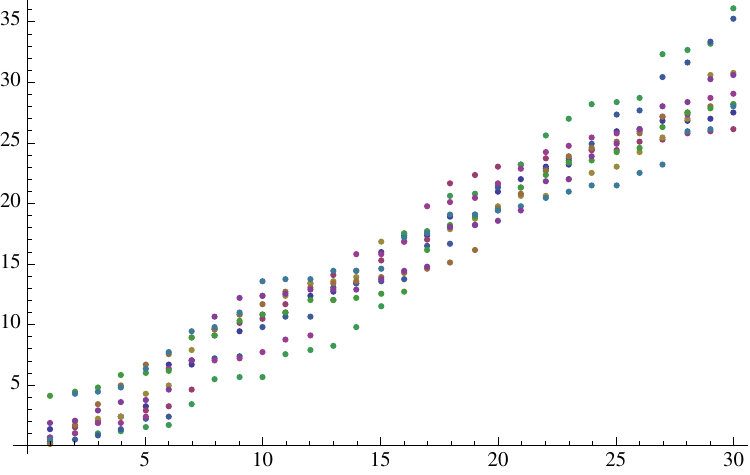}
\caption{Ten random dice with real entries.}
\label{sample-plots-2}

\end{figure}

For another model for random $n$-sided dice we take $n$ random numbers in the unit interval and sort them into increasing order. Then we rescale them, first by dividing by their total and then multiplying by $n(n+1)/2$, so that now the total is the same as for proper $n$-sided dice. These random dice look a lot more like the proper dice but they still have some values greater than $n$. Figure \ref{sample-plots-2} shows a sample of ten of them with $30$ sides. We generated $1000$ triples of these dice with $n=30$ and got $130$ triples with one or more ties. Of the remaining $870$ there were $151$ intransitive sets, giving a ratio of $151/870 \approx 0.174$. There is less intransitivity for these dice than for proper dice but still much more than for the improper dice. Random samples with more sides show the ties decreasing and the proportion of intransitive staying around $0.17$ or $0.18$. We do not have enough evidence to hazard a guess for the limiting value of the proportion. 

Finally, if we consider $n$-sided dice with face numbers from $1$ to $n$ but no restriction on the total, then random samples of three such dice almost never produce ties or intransitive triples. For example, in one run of $1000$ triples of $50$-sided dice there were three triples with a tie and three intransitive triples. 

\section*{One-step dice}

The dice that are closest to the standard die are obtained by moving a single pip from one face to another. That is, the value on one face increases by one and the value on another face decreases by one. If we define the distance between two dice $A$ and $B$ to be 
\[ \sum_i |a_i-b_i|, \]
then these dice are the minimal distance of $2$ from the standard die. We call them ``one-step dice," because they are one step away from the standard die in the graph whose vertices are proper dice and edges between nearest neighbors. 

Let $s(a,b)$ denote the one-step die in which side $a$ goes up by $1$ and side $b$ goes down by $1$. For example, with $n=8$, 
\[  s(2,5) = (1,3,3,4,4,6,7,8)  \mbox{ and }  s(5,2)= (1,1,3,4,6,6,7,8) .\]
In the first, the $2$ changes to $3$ and the $5$ changes to $4$. In the second the $5$ changes to $6$ and the $2$ changes to $1$. The die $s(a,b)$ has a repeated value of $a+1$ and a repeated value of $b-1$, so that $s(a,b)$ has two pairs of repeated values unless $a+1=b-1$, in which case it has one value repeated three times. 

Now the number of one-step dice is much smaller than the number of proper dice. We leave it to the reader to verify that the number of one-step $n$-sided dice is $(n-2)^2$. With such a restricted ensemble of dice, we wondered whether we could understand the prevalence of intransitive sets more completely than for all proper dice. However, for one-step dice ties are common. The one-step dice are not much different from the standard die, and the standard die ties all other proper dice, a fact that we'll need in the next proof. (To see that the standard die ties everyone else, use the representation of proper dice in the heuristic proof of Conjecture 1.)

\begin{proposition}
As $n \rightarrow \infty$ the probability of a tie between a random pair of one-step dice goes to $1$. 
\end{proposition}

\begin{proof}
We consider what happens in the comparison between two dice when we change the value on one face of one die. Suppose that with $A$ we increase the value by one on a single face by replacing $i$ by $i+1$. If $B$ has one face with the value $i$ and one face with value $i+1$, then in the tally of comparisons between all faces of the two dice, there is a net increase of one win for $A$. Similarly, if we reduce one face of $A$ by one, say from $i$ to $i-1$, and if $B$ has a one face with $i$ and once face with $i-1$, then $A$ has a net increase of one loss. We first compare the standard die with $s(c,d)$, which is a tie because the standard die ties every proper die. Now we change the standard die to make it $s(a,b)$. The result is a tie as long as $a$, $a+1$, $b$, and $b-1$ are not repeated values for $s(c,d)$.  Thus, the two dice will tie if  $|a-c|, |a-d|, |b-c|, |b-d| > 2$. As $n$ increases and the values of $a,b,c,d$ are selected randomly, the probability that these inequalities hold approaches $1$.
\end{proof}

So, we have a major difference between one-step dice and all proper dice: as $n$ grows ties become more likely for one-step dice and less likely for proper dice. On the other hand, when we just look at triples of one-step dice in which there are no ties, we see the same behavior as for proper dice: very close to one fourth of the triples are intransitive. For $n=10$ there are $64$ one-step dice and ${64\choose 3 }= 41{,}664$ sets of three dice, of which $8086$ have no ties. There are $2072$ intransitive sets, a proportion of $0.256$. With $n=20$ there are $324$ one-step dice, and ${324 \choose 3 }= 5{,}616{,}324$ sets of size three. We randomly sampled $100{,}000$ sets and found $3664$ with no ties, $907$ of them intransitive, for a proportion of $0.2475$. 

Next we analyze the four scenarios in which one of the modified faces of $s(a,b)$ is close to one of the modified faces of $s(c,d)$ to find out what must hold so that $s(a,b) \succ s(c,d)$. For example, consider the possibility that $a$ and $c$ are close, which means $|a-c| \le 2$. We can assume that $a \ge c$ without loss of generality. If $a=c$, there is a tie. Now consider $a=c+1$. The first die $s(a,b)=s(c+1,b)$ now has among its face values the sequence $c-1, c, c+2, c+2$ while the second die $s(c,d)$ contains the sequence $c-1,c+1,c+1,c+2$. (The rest of the values of the two dice are not relevant.) In the $16$ pairwise comparisons of these, the first die wins seven, loses six and ties three. Therefore, $s(c+1, b)$ beats $s(c,d)$. The other possibility is that $a=c+2$. The die $s(a,b)=s(c+2,b)$ has the face values $c,c+1,c+3,c+3$, while $s(c,d)$ contains $c-1,c+1,c+1,c+2$. Now the $16$ pairwise comparisons result in six wins for each die and four ties, and so $s(c+2,b)$ and $s(c,d)$ are of equal strength. 

By analyzing each of the other three possibilities for $a$ or $b$ interacting with $c$ or $d$, we establish the following lemma.

\begin{lemma} In order for $s(a,b)$ to dominate $s(c,d)$, one or more of the following must hold:
\begin{align*}
  a &= c+1 \\
  d &= a+2 \\
  b &= c \\
  b &= d+1
\end{align*}
\end{lemma}

\begin{proposition}
If $A,B,C$ are randomly chosen one-step dice with no ties among them such that $A\succ B\succ C$, then in the limit as $n \rightarrow \infty$, the two outcomes $A\succ C$ (transitive) and $C \succ A$ (intransitive) are equally likely. Consequently, if three randomly chosen one-step dice have no ties among them, then in the limit as $n \ra \infty$ the probability that they are intransitive approaches $1/4$. 
\end{proposition}

\begin{proof}
With the lemma and the help of a computer program we can estimate the number of solutions for the two alternatives:
\begin{align*} 
    &s(a,b) \succ s(c,d) \succ s(e,f) \prec s(a,b)  \quad \textrm{(transitive)}\\
    &s(a,b) \succ s(c,d) \succ s(e,f) \succ s(a,b)  \quad  \textrm{(intransitive)}
\end{align*}
From the lemma we see that each comparison can occur in four ways. Each alternative requires three comparisons, and so there are potentially $4^3$ scenarios for each. However, some of them are logically impossible;  for example, in the intransitive alternative, the choices $a=c+1$, $c=e+1$, and $e=a+1$ lead to the contradiction $a=a-2$. Now each scenario is represented by a system of three linear equations in the six variables $a,b,c,d,e,f$. Our computer program checks each of the $64$ systems to find those that have positive integer solutions corresponding to one-step dice. The result is that for each alternative $47$ of the $64$  are feasible.  

Because of boundary effects, the scenarios do not have exactly the same number of solutions, but they each have on the order of $n^3$ solutions, since there are three free variables. The boundary effects result in a lower order correction to the dominant $n^3$ term. Therefore, the number of solutions for each alternative is on the order of $47n^{3}$, and so in the limit the two alternatives are equally likely.
\end{proof}

\section*{The big conjectures}

We have seen that intransitive sets of three dice are actually quite common, but what about longer cycles of intransitive dice? Do they even occur? Is there a maximal length? In 2007 Finkelstein and Thorp \cite{FT07} gave an explicit construction of intransitive cycles of arbitrary length. For  example, their construction gives an intransitive cycle of length $5$
\[ A_{1}\succ A_{2} \succ A_{3} \succ A_{4} \succ A_{5} \succ A_{1} \]
with these $15$-sided dice:
\begin{align*}
    A_{1}&=(7,7,7,7,7,7,7,7,7,7,7,7,12,12,12) \\
    A_{2}&=(6,6,6,6,6,6,6,6,6,11,11,11,11,11,11) \\
    A_{3}&=(5,5,5,5,5,5,10,10,10,10,10,10,10,10,10) \\
    A_{4}&=(4,4,4,9,9,9,9,9,9,9,9,9,9,9,9) \\
    A_{5}&=(8,8,8,8,8,8,8,8,8,8,8,8,8,8,8)
\end{align*}
For each odd integer they exhibit an intransitive cycle of that length consisting of dice with three times as many sides. To get a cycle of even length, just construct a cycle of length one greater and delete one of the dice. 

How common are intransitive cycles? With four dice they are quite common. Here are the results from random samples of 1000 sets of four dice having $50, 100, 150, 200$ sides.
\[ \begin{array}{crc}
n & \mbox{ Ties} & \mbox{Intransitive}\\
\hline
50 & .061 & .359 \\
100 & .029 & .365  \\
150 & .023 & .381  \\
200 & .008 & .392 
\end{array}
\]

It definitely looks like the probability of there being any ties goes to $0$, but it's less clear what is happening to the intransitive probability. Before reading further you might want to hazard a guess as to the limiting probability that four random dice form an intransitive cycle.

We have some far-reaching conjectures that go far beyond three or four dice. As consequences we can conjecture the probability (in the limit) that a random set of $k$ dice form an intransitive cycle or that they form a completely transitive set. These conjectures also imply that for proper dice the dominance relation exhibits no bias in favor of transitivity as the number of sides goes to infinity. We consider $k$ random $n$-sided proper dice $A_{1}, A_{2},\ldots,A_{k}$ for a fixed integer $k \ge 2$.

\begin{conjecture}
The probability that there is a tie between any of the $k$ dice goes to $0$ as $n \ra \infty$.
\end{conjecture}

When there are no ties between any of the dice, then there are $2^{\binom{k}{2}}$ outcomes for all the pairwise comparisons among the dice, and each of these outcomes is represented by a tournament graph on $k$ vertices. The vertices are $A_{1},A_{2},\ldots,A_{k}$ and there is an edge from $A_{i}$ to $A_{j}$ if $A_{i} \succ A_{j}$. (A \textbf{tournament graph} is a complete directed graph and is so-called because it represents the results of a round robin tournament.) There are $2^{\binom{k}{2}}$ tournament graphs.
  
\begin{conjecture}
In the limit as $n \ra \infty$ all the tournament graphs with $k$ vertices are equally probable.
\end{conjecture}

Let's apply this conjecture to the case of four dice. There are six comparisons among the pairs of dice and so there are $2^{6}=64$ different tournament graphs. How many of these graphs contain a cycle of length $4$? There are six ways to cyclically arrange the four vertices. Then the remaining two edges can point in either direction. Thus, there are $24$ tournament graphs that contain a $4$-cycle. Therefore, the probability of an intransitive cycle should go to $24/64=3/8$. The experimental results are consistent with the $3/8$ conjecture.

Similar reasoning predicts that the probability of a completely transitive arrangement of four dice has a limit of $3/8$, because there are $4!$ tournament graphs that allow the vertices to be linearly ordered. There are two more symmetry classes of four vertex tournament graphs. In each there is a $3$-cycle with the fourth vertex either dominating or dominated by the vertices in the $3$-cycle. There are $8$ tournament graphs in each of these symmetry classes. Our simulation results are consistent with the prediction that the probability of a completely transitive set is $3/8$ and for the other two classes the probabilities are each $1/8$.

Under the assumption that Conjecture 4 holds, you can predict the probability that a random set of $k$ dice forms an intransitive cycle by finding the number of tournament graphs that contain a cycle through all the vertices, i.e., a Hamiltonian or spanning cycle. Let $C(k)$ be the number of such tournament graphs. The predicted probability is then 
\[   \frac{C(k)}{2^{\binom{k}{2}}} .\]
Basic information about these graphs can be found in the classic book by Moon \cite{Moon68}, where it is shown that having a spanning cycle is equivalent to two other properties: \textbf{strongly connected} or \textbf{irreducible}.  Let $C(k)$ be the number of tournament graphs with $k$ vertices that have a spanning cycle. The $C(k)$ satisfy the equation 
\[ \sum_{j=1}^{k} C(j) 2^{\binom{k-j}{2}}= 2^{\binom{k}{2}} ,\]
and so they can be computed recursively. We have already seen that $C(3)=2$ and $C(4)=24$. Using these values and $C(1)=1$ and $C(2)=0$, we find that $C(5)=544$. Thus, we expect the probability that five random dice are intransitive to approach $544/2^{10}= 17/32 \approx 0.531$ as the number of sides increases. (The sequence $C(k)$ appears in the Online Encyclopedia of Integer Sequences \cite{OEIS-2} as the ``number of strongly connected labeled tournaments.''  )

Wouldn't you guess that the more dice you have the less likely it should be that they are intransitive? But what we are seeing is exactly the opposite. And, in fact, for tournament graphs Moon and Moser proved in 1962 \cite{MoonMoser62} that as $k \ra \infty$ the proportion with spanning cycles goes to $1$. Already for $k=16$ the proportion exceeds $0.999$. 

So we end up with our original beliefs turned on their heads. The dominance relation for proper dice not only fails to be transitive, it is almost as far from transitive as a binary relation can be. We do not know of any other natural example of a binary relation that shows this behavior. Furthermore, our intuition that intransitive dice are rare and that larger sets are even rarer is completely unfounded. They are common for three dice and almost unavoidable as the number of dice grows.

\vspace{4mm}
\noindent \textbf{Acknowledgment} \quad We would like to express our appreciation to Byron Schmuland and Josh Zucker for their helpful discussions.

\end{document}